\documentclass{amsart} %article
\usepackage{amssymb,latexsym}

\newtheorem{thm}{Theorem}[section]
\newtheorem{cor}[thm]{Corollary}
\newtheorem{lem}[thm]{Lemma}

\newtheorem{prop}[thm]{Proposition}
\newtheorem{defn}[thm]{Definition}
\newtheorem{rem}[thm]{Remark}

\title{Moduli space of general connections}

\author{Stanislav Dubrovskiy} %dubr@neu.edu }

\address{Department of Mathematics, Northeastern University, Boston Mass. 02115}

\email{dubr@neu.edu}


\date{July 17, 2010}

\begin{document}

\maketitle

\begin{abstract}
We consider local invariants of  general connections (with torsion).
The group of origin-preserving diffeomorphisms acts on a space of
jets of general connections. Dimensions of moduli spaces of generic connections
are calculated. Poincar\'e series of the geometric structure of connection is
constructed, and shown to be a rational function, confirming the finiteness assertion of Tresse.
\end{abstract}

%%% ----------------------------------------------------------------------

%%% ----------------------------------------------------------------------

\section{Introduction}
A problem of finiteness of functional moduli in various
local differential-geometric settings was discussed by Arnol'd in \cite{A}.
We consider general (not necessarily torsion-free) connections, under the action of smooth coordinate changes.

The structure of the resulting moduli space is reflected in the corresponding Poincar\'e series, explicitly calculated for the generic case (\ref{eq:P'series}). This series turns out to be a rational function (\ref{P:explicit}),
indicating a finite number of invariants. %:
%\begin{equation}
%%\label{eq:P'series}
%=n\sum_{k=1}^{\infty}\;\biggl[\frac{n(n+1)}{2}{n+k-1 \choose n-1}
%-n{n+k+1 \choose n-1}\biggr]t^{k}
%\end{equation}
%$$
%\label{P'ratio}
%p_{\Gamma}(t)=nD_{\Gamma}(\frac{1}{1-t})\ ,
%$$
%where $$D_{\Gamma}=\frac{n(n+1)}{2}{n+t\frac{d}{dt}-1 \choose n-1}
%-n{n+t\frac{d}{dt}+1 \choose n-1}$$
This confirms the finiteness assertion
of Tresse \cite{T} formulated for any ``natural'' differential-geometric structure.

Rationality was earlier confirmed for a symmetric (torsion-free) connection in \cite{D1},
and via a different approach (used also to obtain the result presented herein) in \cite{D2}; and for Fedosov structure (a.k.a. symplectic connection) in \cite{D3}, by the author.

Similar results for Riemannian, K\"{a}hler and hyper-K\"{a}hler structures were
obtained in \cite{Sh1}, and an explicit normal form for
Riemannian structure - in \cite{Sh2}, by Shmelev. Earlier Vershik and Gershkovich investigated
jet asymptotic dimension of moduli spaces of jets of generic distributions at origin in $\mathbb{R}^{n}$
in \cite{VG}, and their normal form in \cite{G}.

I would like to acknowledge support for this work from M.Shubin, and hospitality of the Center of Advanced Studies in Mathematics at Ben Gurion University. I am grateful to C.-L.Terng and A.Vlassov for fruitful discussions.

%------------------------'ere you go--------------------------------------
%we consider coordinate changes on connections, and are in .

\section{Preliminaries and main result}
Let $\mathcal{F}$ and $\mathcal{F}_{k}$ be spaces of germs and $k$-jets respectively 
of smooth connections at the origin in $\mathbb{R}^{n}$.
Two smooth functions
are said to have the same $k$-jet at the origin in $\mathbb{R}^{n}$  if their first $k$ derivatives are equal in some (hence any) local coordinates.

Two connections $\nabla$ and $\tilde{\nabla}$ are said to have the same $k$-jet at 0 if for any smooth vector fields $X$, $Y$, and any smooth function $f$, the functions $\nabla_{X}Y(f)$
and $\tilde{\nabla}_{X}Y(f)$ have the same $k$-jet at $0$. (This is equivalent to Christoffel symbols of
$\nabla$ and $\tilde{\nabla}$ having the same $k$-jet.)\\
%This notion is in fact invariant,
%i.e. if another coordinate system is introduced, Christoffel symbols change, yet
%the equivalence relation holds.
We will frequently denote the connection and its Christoffel symbol with the same letter,
e.g. $\Gamma$, $j^{k}\Gamma$ would stand for its $k$-jet .

There is an action of the group of germs of origin-preserving diffeomorphisms
$G:=\mathrm{Diff}(\mathbb{R}^{n},0)$ on $\mathcal{F}$ and $\mathcal{F}_{k}$.
For $\varphi \in G$, $\nabla$(or $\Gamma$) $\in \mathcal{F}$ and $j^{k}\Gamma \in \mathcal{F}_{k}$:
$$\Gamma \mapsto \varphi^{*}\Gamma\,, \quad j^{k}\Gamma \mapsto j^{k}(\varphi^{*}\Gamma)\,,$$
where $${( \varphi^{*}\nabla )}_{X} Y={\varphi}_{*}^{-1}( \nabla_{  \phantom{|}\atop
{\varphi}_{*}X  } {\varphi}_{*}Y )$$
Let us introduce a filtration of $G$ by normal subgroups:
$$ G=G_{1}\rhd G_{2}\rhd G_{3}\rhd \ldots,$$
where $$G_{k}=\{
\,\varphi \in G\ |\ \varphi(x)=x+(\varphi_{1}(x),\ldots\varphi_{n}(x)),\ \varphi_{i}=O(|x|^{k}),
\,i=1,\ldots,n\,
\}\ .
$$
The subgroup $G_{k}$ acts trivially on $\mathcal{F}_{p}$ for $k\geq p+3$.
It means that the action of $G$ coincides with that of $G/G_{p+3}$ on each $\mathcal{F}_{p}$.
Now $G/G_{p+3}$ is a finite-dimensional Lie group, which we will call $K_{p}$.
%$j^{k}\Gamma\in\mathcal{F}_k$
Denote by $\mathrm{Vect}_{0}(\mathbb{R}^{n})$ the Lie algebra of $C^{\infty}$-vector fields,
vanishing at the origin.
%Define the action of algebra of origin-preserving vector fields
It acts on $\mathcal{F}%_{k}
$ as follows:
%with modified Lie derivative.
\begin{defn}
%\label{d:LD}
For $V\in\mathrm{Vect}_{0}(\mathbb{R}^{n})$
generating a local 1-parameter subgroup $g^{t}$ of $
\mathrm{Diff}(\mathbb{R}^{n},0)$,
\emph{the Lie derivative} of a connection $\nabla$ in the direction $V$ is a (1,2)-tensor:
$$
\mathcal{L}_{V}\nabla%\ (X,Y)
=\left. \frac{d}{dt}\right|_{t=0}
%\,
{g^{t}}^{*}\nabla%\ (X,Y )
$$
\end{defn}
%Before we introduce this action, let us recall how vector field acts %(=Lie differentiates)
%on a connection.
%Connections make up an affine space, with space of (1,2)-tensors as its
%underlying space. Therefore, Lie derivative of a connection is a (1,2)-tensor:
%\begin{lem}
\begin{lem}
\begin{equation}
\label{eq:LD}
(\mathcal{L}_{V}\nabla)(X,Y)=[V,\nabla_{X}Y]-\nabla_{[V,X]}Y-\nabla_{X}[V,Y]
\end{equation}
%\end{lem}
\end{lem}
{\bf Proof}\quad
Below the composition $\circ$ is understood as that of differential operators acting on functions.
$$
(\mathcal{L}_{V}\nabla)(X,Y)=\left. \frac{d}{dt}\right|_{t=0}
g_{*}^{-t}[ \nabla_{  \phantom{|}\atop { g^{t}_{*}X } } g^{t}_{*}Y  ]=
\left. \frac{d}{dt}\right|_{t=0}
\left[\,(g^{t})^{*}\circ [\nabla_{  \phantom{|}\atop { g^{t}_{*}X } } g^{t}_{*}Y ] \circ (g^{-t})^{*}\right]=
$$
$$
\left. \frac{d}{dt}\right|_{t=0} (g^{t})^{*} \circ \nabla_{X}Y +
\nabla_{X}Y \circ \left. \frac{d}{dt}\right|_{t=0} (g^{-t})^{*} +
\nabla_{ \left. \frac{d}{dt}\right|_{t=0}g^{t}_{*}X }Y+
\nabla_{X}\left. \frac{d}{dt}\right|_{t=0}g^{t}_{*}Y =
$$
$$
V\circ \nabla_{X}Y - \nabla_{X}Y \circ V - \nabla_{\left. \frac{d}{dt}\right|_{t=0}g_{*}^{-t}X}Y
 - \nabla_{X}\left. \frac{d}{dt}\right|_{t=0}g_{*}^{-t}Y =
$$
$$
\mathcal{L}_{V}(\nabla_{X}Y)-\nabla_{\mathcal{L}_{V}X}Y - \nabla_{X}(\mathcal{L}_{V}Y)
$$
    \hfill$\Box$\\
This defines the action on germs of connections.

Now we can define the action of $\mathrm{Vect}_{0}(\mathbb{R}^{n})$ on the space of jets $\mathcal{F}_{k}$.\\
For $V\in \mathrm{Vect}_{0}(\mathbb{R}^{n})$:
$$\mathcal{L}_{V}(j^{k}\Gamma)=j^{k}(\mathcal{L}_{V}\Gamma)\ ,$$
where $\Gamma$ on the right
is an arbitrary representative of the $j^{k}\Gamma$ on the left.\\
This is well-defined, since in the coordinate version of (\ref{eq:LD}):
\begin{equation}
\label{eq:LDcoor}
(\mathcal{L}_{V}\Gamma)_{ij}^{l}=
V^{k}\frac{\partial \Gamma_{ij}^{l}}{\partial x^{k}}-
\Gamma_{ij}^{k}\frac{\partial V^{l}}{\partial x^{k}}+
\Gamma_{kj}^{l}\frac{\partial V^{k}}{\partial x^{i}}+
\Gamma_{ik}^{l}\frac{\partial V^{k}}{\partial x^{j}}+
\frac{\partial^{2}V^{l}}{\partial x^{i} \partial x^{j}}
\end{equation}
%\nopagebreak
elements of $k$-th order and less are only coming from $j^{k}\Gamma$, because $V(0)=0$.
Einstein summation convention in (\ref{eq:LDcoor}) above and further on is assumed.\\
Consequently, the action is invariantly defined.
This can also be expressed as commutativity of the following diagram:
$$
\begin{array}{ccccccccccccc}
j^{0}\mathcal{F}& \longleftarrow & \ldots & \longleftarrow & j^{k-1}\mathcal{F} & \stackrel{\pi_{k}}{\longleftarrow} & j^{k}\mathcal{F} &
\longleftarrow & \ldots & \longleftarrow & \mathcal{F} &  &
\\
\downarrow\lefteqn{\mathcal{L}_{V}}& & & & \downarrow\lefteqn{\mathcal{L}_{V}} & & \downarrow\lefteqn{\mathcal{L}_{V}} & &
 & & \downarrow\lefteqn{\mathcal{L}_{V}}& &\!\!\!,
\\
j^{0}\Pi & \longleftarrow & \ldots & \longleftarrow & j^{k-1}\Pi & \stackrel{\pi_{k}}{\longleftarrow} & j^{k}\Pi  &
\longleftarrow & \ldots & \longleftarrow & \Pi & &
\\
\end{array}
$$
%\vskip-0.5cm
where $\pi_{k}$ is the projection from $k$-jets onto $(k-1)$-jets,
$\mathcal{F}$ and $\Pi$
denote the spaces of germs of connections and that of (1,2)-tensors,
respectively, at $0$.

\emph{Poincar\'e series} will encode information about these actions for all $k$.\\
The space $$\mathcal{M}=\mathcal{F}/\mathrm{Diff}(\mathbb{R}^{n},0)$$ of
$\mathrm{Diff}(\mathbb{R}^{n},0)$-orbits on $\mathcal{F}$ is called the \emph{moduli space}
of connections at $0$ on $\mathbb{R}^{n}$.\\ We do not introduce any topology on $\mathcal{M}$.
Similarly, the orbit space
$$\mathcal{M}_{k}=\mathcal{F}_k/\mathrm{Diff}(\mathbb{R}^{n},0)=\mathcal{F}_k/K_{k}$$
is called the moduli space of connection $k$-jets.

The action of $K_k$ is algebraic, a subspace $\mathcal{F}_{k}^{0}\subset\mathcal{F}_{k}$ of points on generic orbits
(those of the largest dimension) is a smooth manifold, open and dense in $\mathcal{F}_{k}$.

A subspace of points on orbits of any other given dimension is a manifold as well,
albeit of a lesser dimension. %, definitely less than a full $dim\mathcal{F}_{k}=dim\mathcal{F}_{k}^{0}$.
We could consider the $G$-quotient for each of those subspaces, and have a moduli space of its own for
each of the orbit types. 

Let $\mathcal{O}_{k}$ denote a generic orbit. Denote by $\mathcal{M}_{k}^{0}$
the moduli space of generic connections:
$$\mathcal{M}_{k}^{0}=\mathcal{F}^0_k/\mathrm{Diff}(\mathbb{R}^{n},0)=\mathcal{F}^0_k/K_{k}\ ,$$
or the generic subspace of the moduli space $\mathcal{M}_k$ . Its
dimension is found as:
\begin{equation}
\label{eq:dim}
\dim\mathcal{M}_{k}^{0}=\dim\mathcal{F}_{k}^{0}-\dim\mathcal{O}_{k}
\end{equation}
%and similarly for other subspaces.
%Remarkably, it turns out
However it is no longer true that the generic moduli (sub)space retains maximal dimension after the quotient is taken, even if we restrict our attention to algebraic Lie group actions. For an explicit counterexample see \cite{MWZ}.\\
%Indeed, for a fixed $x\in\mathcal{F}_{k}$, consider the action map
%$$\varphi_x:G\times\{x\}\rightarrow\mathcal{F}_{k}\,$$
%$$g\mapsto g\cdot x\,,$$
%and its derivative map at identity. Denote by $\mathcal{F}_k^{j}$ the subspace of points on
%orbits of dimension $j$ less than a generic orbit ( which has the maximal dimension ).
%If $x\in\mathcal{F}_k^{j}$, that means that $\rm{rank}\: D\varphi(x)$ is $j$ less than the maximal possible, which
%entails at least $j$ independent conditions annihilating the minors of $D\varphi$.
%Hence codim$\:\mathcal{F}_k^{j} \geq j$, so even though orbit dimension may drop by $j$, it must be at least matched
%by the drop in subspace dimension, compared to generic case, and generic moduli space retains maximal dimension.
Thus we define $$\dim\mathcal{M}_{k}=\dim\mathcal{M}_{k}^{0}%=dim \mathcal{F}_{k}-dim \mathcal{O}_{k}\; ,
$$ for mere simplicity of notation.
One more piece of notation:
\begin{displaymath}
a_{k}=\left\{ \begin{array}{ll} \dim \mathcal{M}_{k}\,, & k=0 \\
                                \dim \mathcal{M}_{k}- \dim \mathcal{M}_{k-1}\,, & k\geq 1
\end{array} \right.
\end{displaymath}
and we can introduce our main object of interest.

\begin{defn}
The formal power series
$$p_{\Gamma}(t)=\sum_{k=0}^{\infty}a_{k}t^{k}$$
is called the \emph{Poincar\'e series} for the moduli space $\mathcal{M}$.
\end{defn}
Our main result is the following theorem.
\newpage
\begin{thm}
\label{thm:series}
Poncar\'e series coefficients $a_{k}=a(k)$ are polynomial in $k$, and the series has the form:\\\\
$p_{\Gamma}(t)=%(\delta_{1}^{n}+\delta_{2}^{n}-n^2)t-\delta_{2}^{n}t^2+
{\displaystyle n\sum_{k=1}^{\infty}\;\biggl[n^2{n+k-1 \choose n-1}
-{n+k+1 \choose
n-1}\biggr]t^{k}+\frac{n^2(n-3)}{2}+\delta_{1}^{n}+2\delta_{2}^{n}(1-t)\,.}
$\\\\
$\hspace*{22em}(\ \delta\textrm{\emph{ is the Kronecker symbol }}).$\\\\
It represents a rational function.
\end{thm}
\begin{rem}
This complies with Tresse' assertion that algebras of ``natural"
differential-geometric structures are finitely-generated.
\end{rem}
{\bf Proof}\ \ of this theorem is relegated to section \ref{sec:Proof}.\\\\
To explain significance of rationality of Poincar\'e series we make the following
\begin{rem}
\label{rem-finiteness-poles-at-1}
If a geometric structure is
described by a finite number of functional moduli, then its
Poncar\'e series is rational. In particular, if there are $m$
functional invariants in $n$ variables, then
$$p(t)=\frac{m}{(1-t)^n}$$
\end{rem}
Indeed, dimension of the moduli spaces of $k$-jets is just the number
of monomials up to the order $k$ in the formal power series of the
$m$ given invariants:$$\dim \mathcal{M}_{k}=m{n+k \choose n}$$ For
more details and slightly more general formulation see {\bf
Theorem 2.1} in \cite{Sh2}.
%------------------------------------------------------------------------------------------------------
\section{Stabilizer of a generic $k$-jet}
%%% ----------------------------------------------------------------------
The main challenge in calculating a Poincar\'e series is finding
the stabilizer of a generic $k$-jet. Here we do it in fixed
normal coordinates. The suggestion to use special coordinates is due to Vlassov, \cite{V}.

Let $\nabla$ be a connection given in local coordinates around origin
by its Christoffel symbols $\Gamma^{i}_{jk}$.

\begin{defn} A geodesic curve is the solution of
\begin{equation}
\label{n:def:geo}
\frac{d^{2}x^i}{d t^{2}}+\Gamma^{i}_{jk}\frac{dx^j}{d t}\frac{dx^k}{d t}=0\ ,\quad x(0)=0\ ,\quad x'(0)=v\ .
\end{equation}
\end{defn}
\begin{defn}
A coordinate system $x$ is called an \textit{affine normal coordinate system} when
solutions of \eqref{n:def:geo} are linear in the canonical parameter $t$:
\begin{equation}
\label{n:eq:def-norm-geo}
x^i=a^it\ .
\end{equation}
\end{defn}
These are exactly the coordinates induced from the tangent space by the
exponential map of $\Gamma$
$$
\textrm{exp}_0(v)=x(1)
\ ,$$
where $v\in U$, a small neighborhood of 0 in $T_0 \mathbb{R}^n$, and $x(t)$ is
the geodesic defined by \eqref{n:def:geo}.
%%%%%%%%%%%%%%%%%%%%   Lemma 4.3   %%%%%%%%%%%%%%%%%%%%%%%%%%%%%%%
\begin{lem}
\label{n:lem:norm}
Coordinates $x$ on $\mathbb{R}^n$ are normal if and only if:
\begin{equation}
\label{n:eq:normal}
\Gamma^{i}_{jk}(x)x^jx^k\equiv0\ .
\end{equation}
\end{lem}
%%%%%%%%%%%%%%%%    C O R O L L A R Y    4 . 4   %%%%%%%%%%%%%%%%%%%%%%%%%%%
\begin{cor}
\label{n:cor:norm}
Let $x$ be a normal coordinate system. We have:
\begin{equation}
\label{n:eq:gamma(0)=0}
\Gamma^{i}_{jk}(0)+\Gamma^{i}_{kj}(0)=0
\end{equation}
\begin{equation}
\label{n:eq:circus}
\left. \left(
\frac{\partial \Gamma^i_{jk} }{ \partial x^l }+\frac{\partial \Gamma^i_{kj} }{ \partial x^l }+
\frac{\partial \Gamma^i_{lj} }{ \partial x^k }+\frac{\partial \Gamma^i_{jl} }{ \partial x^k }+
\frac{\partial \Gamma^i_{kl} }{ \partial x^j }+\frac{\partial \Gamma^i_{lk} }{ \partial x^j }
\right) \right|_{x=0}=0
\end{equation}
for distinct $j, k, l$ (i.e. terms with repeating lower indexes, as $\begin{displaystyle}\frac{\partial \Gamma^i_{jj} }{ \partial x^k }\end{displaystyle}$, %are to 
appear only once).

In general for $r=1, 2, \ldots$,
\begin{equation}
\label{n:eq:Bigcircus}
S_{j,k\hookrightarrow j,k,\alpha_1 \ldots \alpha_r}
\left(
{
\left.
\frac{ \partial^r \Gamma^i_{jk} }{ \partial x^{\alpha_1} \ldots \partial x^{\alpha_r} }
\right|
}_{x=0}
\right)=0\ ,
\end{equation}\\
where $S_{j,k\hookrightarrow j,k,\alpha_1 \ldots \alpha_r}$ is the
sum of all terms obtained from the term inside the parentheses by
replacing the pair ($j$, $k$) by any ordered pair from the set
$\{j,k,\alpha_1 \ldots \alpha_r\}$.
\end{cor}
In fact, the converse holds as well, cf. \cite[VI.41]{Th} and \cite[Prop.4.2]{GRS}.\\
Namely, let \{$\Gamma^i_{jk\alpha_1 \ldots \alpha_r}|r=1, 2, \ldots $\} be a set of numbers 
(or tensors given at a point),
symmetric in $\{\alpha_1 \ldots \alpha_r\}$, satisfying (analogue of) \eqref{n:eq:gamma(0)=0} and
$$S_{j,k\hookrightarrow j,k,\alpha_1 \ldots \alpha_r}
\left(
\Gamma^i_{jk\alpha_1 \ldots \alpha_r}
\right)=0$$ as in \eqref{n:eq:Bigcircus}, and such that the power series
$$
\sum_{r=1}^\infty \frac{1}{r!}\Gamma^i_{jk\alpha_1 \ldots \alpha_r}x^{\alpha_1}\ldots x^{\alpha_r}
$$
converges near 0.

Then this series defines an object $\Gamma^i_{jk}(x)$ satisfying \eqref{n:eq:normal}
in given coordinates $x$. Hence $\Gamma^i_{jk}(x)$ is a set of Christoffel symbols of a general connection
in given local coordinates $x$, such that the system of coordinates $x$ is normal for it.
We say that Cor.\ref{n:cor:norm}\hspace{.5em}gives \textit{a complete set of identities for Christoffel symbols of a connection in normal coordinates}.

\noindent For more details on normal coordinates, including proofs of the Lemma \ref{n:lem:norm} and its corollary,
consult \cite{Th}, \cite{Veb} or \cite[Sec.4]{GRS}.\\

Let us now turn to the problem of finding the stabilizer of a given connection $\Gamma$.
From this point on through the rest of the section we assume that local coordinates $x$ are normal coordinates associated to $\Gamma$. 

In these coordinates let us introduce grading in homogeneous components on $\Gamma$:
$$\Gamma=\Gamma_{0}+\Gamma_{1}+\ldots\ ,$$
and on diffeomorphism generating $V\in\mathrm{Vect}_{0}(\mathbb{R}^{n})$:
$$V=V_{1}+V_{2}+\ldots$$
\begin{center} ($V_{0}=0$, since $V$ preserves the origin).
\end{center}
A diffeomorphism preserving $\Gamma$ must preserve its set of geodesics,
which in this coordinate system are lines \eqref{n:eq:def-norm-geo}, all parametrized by the canonical parameter $t$. Such diffeomorphisms are linear maps, e.g.:
$$V=V_1\ , \quad V_2=V_3=\ldots=0\ .$$
Then the stabilizer condition $\mathcal{L}_{V}\Gamma =  0$ considered as an equation on $V$ becomes:\\
\begin{equation}
\label{n:sys:SYS}
\left\{
%{\setlength\arraycolsep{2pt}
\begin{array}{l}
\mathcal{L}_{V}\Gamma_{0} =  0 \\
\mathcal{L}_{V}\Gamma_{1} =  0\\
\qquad \vdots  \\
\mathcal{L}_{V}\Gamma_{k} =  0\;.
\end{array}
\right.
\end{equation}\\
We need to find all $V$ solving this system for a generic
$\Gamma$. Write 
\begin{equation}
\label{notation}
\gamma^{l}_{ij}:=\Gamma^{l}_{ij}(0)\ , \quad V^{k}=\sum_{s=1}^{n}b_{s}^{k}x^{s}\ . 
\end{equation}\\
In this notation the first equation of \eqref{n:sys:SYS} is a linear system on $b_{s}^{k}$:
\begin{equation}
\label{ijl-sys}
-\gamma^{k}_{ij}b_{k}^{l}+\gamma^{l}_{kj}b_{i}^{k}+\gamma^{l}_{ik}b_{j}^{k}=0\ ,
\end{equation}
indexed by $\{(ijl), i<j\}$, it is sufficient to consider $i<j$ since exchanging them leads to the same equation.
The coefficients $\gamma$ are subject only to:
\begin{equation}
\label{asym-cond}
\gamma^{l}_{ij}=-\gamma^{l}_{ji}
\end{equation}\\
due to \eqref{n:eq:gamma(0)=0}.

We will presently show that for $n\geq
3$ this system is non-degenerate in general position. Thus the stabilizer is trivial (for any order jet). Exceptional dimensions 1 and 2 are
considered in the next section. 

It suffices to present one such connection (or just its constant part $\{\gamma^{l}_{ij}\}$) which makes the system non-degenerate.
Namely, set the coefficients to:
\begin{equation}
\label{gamma:conditions}
\gamma^{\alpha}_{\alpha\beta}=
\left[
\begin{array}{l}
\ \;\,1\,,\alpha<\beta \\
-1\,,\alpha>\beta \\
\end{array}
\right.
\end{equation}
and zero otherwise (when the upper index fails to match either of the lower ones).
 Then \eqref{ijl-sys} splits into the following non-trivial cases according to the index $(ijl)$:\\\\
 $l<i<j$\hspace{11ex}$-2b^l_j=0$,\\
 $i<j<l$\hspace{11ex}$-2b^l_i=0$,\\ 
 implying $b^\alpha_\beta=0$ for $\alpha\neq\beta$, and finally\\\\
$l=i<j$\hspace{13ex}$b^j_j=0$,\\
$i<j=l$\hspace{13ex}$b^i_i=0$,\\\\
finishing the proof in dimensions 3 and higher (availability of three distinct indexes is essential).

%--------------------------------------------------------------------------------------------------
\section{Exceptions: the stabilizer in low dimensions}
We start with the case $n=1$. In this case we have only one function in one variable: $$\Gamma^1_{11}(x)\equiv 0\,,$$
because of \eqref{n:eq:normal}.
This means the stabilizer conditions \eqref{n:sys:SYS} are all empty, and the stabilizer has the maximal possible dimension 1.\\

In the case $n=2$, the stabilizer of the 0-jet is determined by \eqref{ijl-sys}.
With only two distinct indexes available it makes for a system of two equations in four variables:\\
%\vspace{0.9ex}
\hspace*{31.7ex}$b_{1}^{1}$\hspace{5ex}$b_{2}^{1}$\hspace{5.25ex}$b_{1}^{2}$\hspace{3.77ex}$b_{2}^{2}$\\
\vspace*{-1ex}
\\$(122)$\\
$(121)$
\vspace*{-5.88ex}
$$
\left(
\begin{array}{ccccc}
%&&&&\\
& \gamma_{12}^{2} & & -\gamma_{12}^{1} &  \\
& & -\gamma_{12}^{2} & & \gamma_{12}^{1}\\
\end{array}\ 
\right)
$$\\
It is of full rank in general position, since the coefficients are only subject to \eqref{asym-cond}, making the stabilizer of the 0-jet 2-dimensional.\\

For the 1-jet, in addition to the above, further restrictions on $V$ are imposed by the second equation in \eqref{n:sys:SYS}.
Denoting $\begin{displaystyle}\frac{\partial \Gamma^i_{jk} }{ \partial x^l }(0)\end{displaystyle}$ by $\Gamma^i_{jk,l}$ , 
this equation produces the following system:
\begin{equation}
\label{ijls-sys}
\Gamma^{l}_{ij,k}b_{s}^{k}-\Gamma^{k}_{ij,s}b_{k}^{l}+\Gamma^{l}_{kj,s}b_{i}^{k}+\Gamma^{l}_{ik,s}b_{j}^{k}=0\ .
\end{equation}
There are 16 equations indexed by arbitrary 4-tuples $(ijls)$. The coefficients satisfy:
$$
\hspace*{-2.3ex}\Gamma^{i}_{11,1}=\Gamma^{i}_{22,2}=0
$$\\
\vspace*{-6.25ex}
\begin{equation}
\label{jet 2-conditions on gamma}
\hspace*{16.21ex}\Gamma^{i}_{12,1}+\Gamma^{i}_{21,1}+\Gamma^{i}_{11,2}=0\qquad\qquad\hspace*{7ex} i=1,2
\end{equation}
$$\Gamma^{i}_{12,2}+\Gamma^{i}_{21,2}+\Gamma^{i}_{22,1}=0\ ,$$
all consequences of \eqref{n:eq:circus}. The twelve nontrivial equations have the matrix below:\\
\vspace*{-1ex}\\
%%%%%% THE M A T R I X %%%%%%%%%
\hspace*{23.9ex}$b_{1}^{1}$\hspace{12.39ex}$b_{2}^{1}$\hspace{15ex}$b_{1}^{2}$\hspace{11.5ex}$b_{2}^{2}$\\
%\vspace*{1ex}
\\
$(1112)$\\
$(1122)$\\
$(1211)$\\
$(1212)$\\
$(1221)$\\
$(1222)$\\
$(2111)$\\
$(2112)$\\
$(2121)$\\
$(2122)$\\
$(2211)$\\
$(2221)$\\
\vspace*{-37ex}
$$
\hspace*{11ex}
\left(
\begin{array}{cccc}
\Gamma_{11,2}^{1} & -\Gamma_{11,2}^{2} & -\Gamma_{22,1}^{1} & \Gamma_{11,2}^{1}\\
2\Gamma^{2}_{11,2}  &0 & -(\Gamma_{11,2}^{1}+\Gamma^{2}_{22,1})&0\\
\Gamma^{1}_{12,1} & -\Gamma^{2}_{12,1} & -\Gamma^{1}_{21,2} & \Gamma^{1}_{12,1}\\
0& -(\Gamma_{21,1}^{1}+\Gamma^{2}_{12,2}) & 0& 2\Gamma_{12,2}^{1}\\
2\Gamma^{2}_{12,1}&0 & -(\Gamma_{21,2}^{2}+\Gamma^{1}_{12,1}) &0 \\
\Gamma^{2}_{12,2} & -\Gamma^{2}_{21,1}& -\Gamma^{1}_{12,2}& \Gamma^{2}_{12,2}\\
\Gamma^{1}_{21,1}& -\Gamma^{2}_{21,1}& -\Gamma^{1}_{12,2}& \Gamma^{1}_{21,1}\\
0& -(\Gamma_{12,1}^{1}+\Gamma^{2}_{21,2})& 0& 2\Gamma^{1}_{21,2}\\
2\Gamma^{2}_{21,1}& 0& -(\Gamma_{12,2}^{2}+\Gamma^{1}_{21,1}) & 0\\
\Gamma^{2}_{21,2}& -\Gamma^{2}_{12,1}& -\Gamma^{1}_{21,2}& \Gamma^{2}_{21,2}\\
0& -(\Gamma_{11,2}^{1}+\Gamma^{2}_{22,1})&0 &2\Gamma^{1}_{22,1} \\
\Gamma^{2}_{22,1}& -\Gamma^{2}_{11,2}& -\Gamma^{1}_{22,1}&\Gamma^{2}_{22,1} \\
\end{array}\ 
\hspace{-1ex}\right)
$$\\
Setting the coefficients $\Gamma^{i}_{12,j}=2$, $\Gamma^{i}_{21,j}=\Gamma^{i}_{kk,j}=-1$ ($i$ and $j$ arbitrary, $k\neq j$) produces a non-degenerate system (e.g. the top four equation minor is non-trivial). 
Thus the only solution is zero and the stabilizer is trivial for jets of order one (and higher) in dimension two.
To summarize:
\begin{prop}\label{stab}
The stabilizer of a k-jet of a generic connections:\\
for $n=1$ is 1-dimensional for any $k$ ;\\
for $n=2$ is 2-dimensional for $k=0$ and trivial for $k\geq 1$ ;\\
for $n\geq 3$ is trivial for any $k$.\\
\end{prop}
%%%%%%%%%%%%%%%%%%% Poincare Series %%%%%%%%%%%%%%%%%%%%%%%%%%%%%
\section{Poincar\'e series}
\label{sec:Proof}
We use the Proposition \ref{stab} above to find the dimension of a generic orbit:\\\\
$
\dim \mathcal{O}_{k}( \Gamma ) = \dim( T_{\mathrm{id}}(K_{k}/G_{\Gamma}))$\\\\
$\hspace*{5em}=\dim(\{ V | V={V}_{1}+{V}_{2}+\ldots+
{V}_{k+1}+{V}_{k+2} \}) -\delta_{1}^{n}-2\delta_{2}^{n}\delta_{0}^{k}\;,
$\\\\
where $\delta$ is a Kronecker symbol, taking care of non-trivial stabilizers for various $k$ and $n$ ;
${V_{i}}$ is an n-component vector,
each component a homogeneous polynomial of degree $i$ in $(x^{1},...,x^{n})$.
So we have:
$$
\dim\mathcal{O}_{k}=
{\displaystyle n\sum_{m=1}^{k+2} {n+m-1 \choose n-1} }-\delta_{1}^{n} \qquad \textrm{for }k\geq1 \;, $$
$$
\dim \mathcal{O}_{0}=
n\sum_{m=1}^{2} {n+m-1 \choose n-1} -\delta_{1}^{n}-2\delta_{2}^{n}=\frac{n^{2}(n+3)}{2}-\delta_{1}^{n}-2\delta_{2}^{n}\ .$$
Dimension of the moduli space of connection $k$-jets $\mathcal{M}_{k}$ is:
$$\dim \mathcal{M}_{k}=\dim \mathcal{F}_{k}-\dim \mathcal{O}_{k}\;,$$
where $\mathcal{F}_{k}$ is the space of connection k-jets.\\
For $k=0$:
$$\dim \mathcal{M}_{0}=\dim \mathcal{F}_{0}-\dim \mathcal{O}_{0}=\frac{n^{2}(n-3)}{2}+\delta_{1}^{n}+2\delta_{2}^{n}\ .$$
For $k\geq1$ :
$$\dim \mathcal{M}_{k}=\dim \mathcal{F}_{k}-\dim \mathcal{O}_{k}$$
$$=n^3\sum_{m=0}^{k}{n+m-1 \choose n-1}-n\sum_{m=1}^{k+2}{n+m-1 \choose n-1}+
\delta_{1}^{n}\ .$$
%%%%%%%%%%%%%%%%%%%%%%%%%%%%%%%%%%%%%%%%%%%%%%%%%%%%%%%%%%%%
The Poincar\'e series is:
\begin{equation}
\label{eq:P'series}
p_{\Gamma}(t)=\dim \mathcal{M}_{0}+\sum_{k=1}^{\infty}(\dim \mathcal{M}_{k}-\dim \mathcal{M}_{k-1})t^{k}
\end{equation}
$$\hspace*{-7ex}
=n\sum_{k=1}^{\infty}\;\biggl[n^2{n+k-1 \choose n-1}
-{n+k+1 \choose
n-1}\biggr]t^{k}+\frac{n^2(n-3)}{2}+\delta_{1}^{n}+2\delta_{2}^{n}(1-t)\,.
$$\\
Note that since $\dim \mathcal{M}_{0}$ is exceptional, the linear term in the series had to be corrected by $-2\delta_{2}^{n}t$.
In dimension one $p_\Gamma(t)\equiv0$.
Simplifying the above, we obtain the following\\\\
{\bf Fact}
\begin{em}
The Poncar\'e series $p_\Gamma(t)$ is a rational function.
Namely,
\begin{equation}
\label{P'ratio}
p_{\Gamma}(t)=
\delta_{1}^{n}+2\delta_{2}^{n}(1-t)-n^2+
nD_{\Gamma}\left(\frac{1}{1-t}\right)\ ,
\end{equation}
where $D_{\Gamma}$ is a differential operator of order $n-1$ :
$$D_{\Gamma}=n^2{n+t\frac{d}{dt}-1 \choose n-1}
-{n+t\frac{d}{dt}+1 \choose n-1}\ ,$$ with $${n+t\frac{d}{dt}-1
\choose
n-1}=\frac{1}{(n-1)!}(t\frac{d}{dt}+1)\ldots(t\frac{d}{dt}+n-1)\
,$$
$${n+t\frac{d}{dt}+1 \choose n-1}=\frac{1}{(n-1)!}(t\frac{d}{dt}+3)\ldots(t\frac{d}{dt}+n+1) \ .$$
\end{em}
Indeed, denote
$$\varphi_{m}(t)=\sum_{k=0}^{\infty}k^{m}t^{k}\ ,\qquad m\in\mathbb{Z}_{+}\ ,$$
then
$$\varphi_{m}(t)=\sum_{k=0}^{\infty}k^{m-1}kt^{k-1}t=t\left(\sum_{k=0}^{\infty}k^{m-1}t^{k}\right)'
=\left(t\frac{d}{dt}\right)\varphi_{m-1}(t)\quad \mathrm{for}\ m\in\mathbb{N}\,.$$
Thus
$$\varphi_{m}(t)=\left(t\frac{d}{dt}\right)^{m}\varphi_{0}(t)=
\left(t\frac{d}{dt}\right)^{m}\left(\frac{1}{1-t}\right)\,.$$
Hence,
$$\sum_{k=0}^{\infty}\;\biggl[n^2{n+k-1 \choose n-1}-{n+k+1 \choose n-1}\biggr]t^{k}$$
$$=\biggl[n^2{n+t\frac{d}{dt}-1 \choose n-1}-{n+t\frac{d}{dt}+1 \choose n-1}\biggr]
\left(\frac{1}{1-t}\right)\,.$$\\
We have to account for the ``extra''
$0^{\textrm{th}}$ term, since the Poincar\'{e} series does not have it:
$$a_{0}=n^2-{n+1 \choose
n-1}=\frac{n(n-1)}{2}\,.$$ 
So
$$p_{\Gamma}(t)=\frac{n^2(n-3)}{2}+\delta_{1}^{n}+2\delta_{2}^{n}(1-t)+
n\sum_{k=0}^{\infty}\;\biggl[n^2{n+k-1 \choose
n-1}-{n+k+1 \choose n-1}\biggr]t^{k}-na_0$$
$$=\delta_{1}^{n}+2\delta_{2}^{n}(1-t)-n^2+
nD_{\Gamma}\left(\frac{1}{1-t}\right) \hspace*{9ex}\Box$$ Now we
shall derive the explicit formula for this rational function.
\begin{lem}
For $N\geq2$:
$$
{N+t\frac{d}{dt}-1\choose N-1}\frac{1}{1-t}=\frac{1}{(1-t)^N}\ ,
$$
$$
{N+t\frac{d}{dt}+1\choose N-1}\frac{1}{1-t}=
\frac{1}{(1-t)^N}+\frac{2}{(1-t)^{N-1}}+\frac{3}{(1-t)^{N-2}}+\ldots
$$
$$\hspace*{16ex}+\frac{k}{(1-t)^{N-k+1}}+
\ldots+\frac{N}{(1-t)}
$$
\end{lem}
{\bf Proof}\quad is by induction on $N$. \hspace{39ex} $\Box$
\\\\
This together with the fact that $p_\Gamma(t)\equiv0$ for $n=1$ implies the following in all dimensions.
\begin{cor}
\begin{equation}
\label{P:explicit}
p_{\Gamma}(t)=\delta_{1}^{n}+2\delta_{2}^{n}(1-t)-n^2+n\left(
\frac{n^2-1}{(1-t)^{n}}-
\frac{2}{(1-t)^{n-1}}-\ldots\frac{n}{(1-t)} \right)\,.
\end{equation}
\end{cor}
This formula shows that $p_{\Gamma}(t)$ is not an arbitrary rational
function, but one of the form required by the Tresse finiteness
claim, with poles exclusively at $t=1$,\\cf. Remark
\ref{rem-finiteness-poles-at-1}.


\begin{thebibliography}{99}

\bibitem[A]{A} Arnol'd, V.I. Mathematical Problems in Classical Physics. \textit{Trends and Perspectives in Applied Mathematics,}
1--20, Appl. Math. Sci., 100, \textit{Springer, New York,} 1994.

\bibitem[D1]{D1} Dubrovskiy, S. Moduli space of symmetric connections, \emph{Zap. Nauchn. Sem. S.-Peterburg. Otdel. Mat. Inst. Steklov. (POMI)}  \textbf{292} (2002),  Teor. Predst. Din. Sist. Komb. i Algoritm. Metody. 7, 22--39, 177;  \emph{translation in  J. Math. Sci. (N. Y.)}  \textbf{126}  (2005),  no. 2, 1053--1063.

\bibitem[D2]{D2} Dubrovskiy, S. Differential invariants of geometric structures.  \textit{The COE Seminar on\\ Mathematical Sciences 2004},  11--36, Sem. Math. Sci., 31, \textit{Keio Univ., Yokohama,} 2004.

\bibitem[D3]{D3} Dubrovskiy, S. Moduli space of Fedosov structures, \emph{Ann. Global Anal. Geom.}  27  (2005),  \\no. 3, 273--297.

\bibitem[G]{G} Gershkovich, V. Ya. On normal form of distribution jets. \textit{Topology and geometry-Rohlin Seminar,} 77--98, Lecture Notes in Math., 1346, \emph{Springer, Berlin,} 1988.

\bibitem[GRS]{GRS} Gelfand I.M., Retakh, V., Shubin, M.A, Fedosov Manifolds,\\
\textit{Advances in Mathematics}, 1998, vol. 136, pp.104-140.

\bibitem[MWZ]{MWZ} P.Magyar, J.Weyman, A.Zelevinsky, Multiple flag varieties of finite type.
\textit{Adv. Math.} 141 (1999), no. 1, 97--118.

\bibitem[Sh1]{Sh1} Shmelev, A.S. On Differential Invariants of Some Differential-Geometric Structures,
\textit{Proceedings of the Steklov Institute of Mathematics}, 1995, vol.209, pp.203-234.

\bibitem[Sh2]{Sh2} Shmelev, A.S. Functional moduli of germs of Riemannian metrics (Russian), \textit{Funktsional. Anal. i Prilozhen.}, 31 (1997), no. 2, pp.58--66, 96; translation in \textit{Funct. Anal. Appl.} 31 (1997), no. 2, pp.119--125

\bibitem[T]{T} Tresse, A. Sur les Invariants Diff\'erentiels des Groupes Continus des Transformations,\\ \textit{Acta Mathematica}, 1894, vol.18, 1-88.

\bibitem[Th]{Th} T.Y.Thomas, ``The Differential Invariants of Generalized Spaces'',\\
Cambridge Univ. Press, London, 1934.

\bibitem[Veb]{Veb} O.Veblen, ``Invariants of uadratic Differential Forms'',
Cambridge Tracts in Mathematics and Mathematival Physics, Vol. 24, Cambridge Univ. Press,\\
Cambridge, UK, 1927 [Reprinted 1952].

\bibitem[VG]{VG} Vershik, A. M.; Gershkovich, V. Ya. Estimation of the functional dimension of the orbit space of germs of distributions in general position. (Russian)
\textit{Mat. Zametki} \textbf{44} (1988), no. 5, 596-603, 700; \textit{translation in Math. Notes} \textbf{44} (1988), no. 5-6, 806-810 (1989).

\bibitem[V]{V} A.T.Vlassov, private communication.
\end{thebibliography}
\end{document}